\RequirePackage{fix-cm}
\documentclass[smallextended]{svjour3}       
\smartqed
\usepackage{graphicx}

\usepackage{amsmath,amscd,amsfonts,amssymb,enumerate, mathrsfs}
\usepackage{color}
\usepackage[colorlinks]{hyperref}

\usepackage{floatrow}
\usepackage{multirow}
\usepackage[table]{xcolor}

\numberwithin{equation}{section}

\numberwithin{equation}{section}

\DeclareMathOperator{\RE}{Re}

\begin{document}
	
	\title{Geometric Properties of Functions Containing Derivatives of Bessel Functions}
	\titlerunning{Geometric Properties of $N_\nu(z)$}       
	
	\author{Sercan Kaz\i mo\u glu \and Kamaljeet Gangania$^*$}
	
	\authorrunning{Sercan Kaz\i mo\u glu and Kamaljeet Gangania} 
	
	\institute{Sercan Kaz\i mo\u glu \at
		\email{srcnkzmglu@gmail.com}  \\
		\emph{ Department of Mathematics, Faculty of Science and Literature, Kafkas University, 36100, Kars-Turkey}\\
		Kamaljeet Gangania \at
		\email{gangania.m1991@gmail.com}  \\
		\emph{ Department of Applied Mathematics, Delhi Technological University, Delhi--110042, India}\\
		${}^*$Corresponding author	    
	}
	
	\date{Received: date / Accepted: date}

	\maketitle
	
	\begin{abstract}
		In this paper our aim is to find the radii of $\gamma$-Spirallike of order $\alpha$ and convex $\gamma$-Spirallike of order $\alpha$ for three different kinds of normalizations of the function $N_\nu(z)=az^2J_\nu^{\prime\prime}(z)+bzJ_\nu^{\prime}(z)+cJ_\nu(z),$ where $J_\nu(z)$ is the Bessel function of the first kind of order $\nu.$ Moreover, the $\mathcal{S}^*\left(\varphi\right)-$radii and $\mathcal{C}\left(\varphi\right)-$radii of these normalized functions are investigated. The
		tables are created and visual verification with graphs are made by giving special values
		to the real numbers $a,$ $b$ and $c$ in the obtained results. 
		
	\end{abstract}	
	\section{Introduction}
	Let $\mathcal{A}$ be the class of analytic functions normalized by the condition $f(0)=0=f'(0)-1$ in the unit disk $\mathbb{D}:=\mathbb{D}_{1}$, where   $\mathbb{D}_{r}:=\{z\in \mathbb{C}: |z|<r \}$.  We say that a function $f\in \mathcal{A}$ is $\gamma$-Spirallike of order $\alpha$ if and only if
	\begin{equation*}
	\RE\left(e^{-i\gamma} \frac{zf'(z)}{f(z)} \right) > \alpha \cos{\gamma},
	\end{equation*}
	where $\gamma \in \left(-\frac{\pi}{2},\frac{\pi}{2}\right)$ and $0\leq\alpha<1$. The class of such functions, we denote by $\mathcal{S}_{p}^{\gamma}(\alpha)$. In view of well known  Alexender's relation, let $\mathcal{CS}_{p}^{\gamma}(\alpha)$ be the class of convex $\gamma$-spirallike functions of order $\alpha$, which is defined below
	\begin{equation*}
	\RE\left(e^{-i\gamma} \left(1+\frac{zf''(z)}{f'(z)}\right) \right) > \alpha \cos{\gamma}.
	\end{equation*}	
	 Spacek~\cite{spacek-1933} introduced and studied the class $\mathcal{S}_{p}^{\gamma}(0)$. Each function in $\mathcal{CS}_{p}^{\gamma}(\alpha)$ is univalent in $\mathbb{D}$, but they do not necessarily be starlike. Further, it is worth to mention that for general values of $\gamma (|\gamma|<\pi/2)$, a function in $\mathcal{CS}_{p}^{\gamma}(0)$ need not be univalent in $\mathbb{D}$. For example: $f(z)=i(1-z)^i-i \in \mathcal{CS}_{p}^{\pi/4}(0)$, but not univalent. Indeed, $f \in \mathcal{CS}_{p}^{\gamma}(0)$ is univalent if $0<\cos\gamma<1/2$, see Robertson~\cite{Robertson-1969} and Pfaltzgraff~\cite{Pfaltzgraff-1975}.  Note that for $\gamma=0$, the classes $\mathcal{S}_{p}^{\gamma}(\alpha)$ and $\mathcal{CS}_{p}^{\gamma}(\alpha)$ reduce to the classes of starlike and convex functions of order $\alpha$, given by
	\begin{equation*}
	\RE\left(\frac{zf'(z)}{f(z)} \right) > \alpha \quad \text{and} \quad
	\RE\left(1+\frac{zf''(z)}{f'(z)}\right) > \alpha,
	\end{equation*}
	which we denote by $\mathcal{S}^*(\alpha)$ and $\mathcal{C}(\alpha)$, respectively.
	
	In the recent past, connections between the special functions and their geometrical properties have been established in terms of radius problems \cite{abo-2018,bdoy-2016,Baricz,Ba3,b-praj-2020,btk-2018,Bricz-Rama,ErhanDenij2020}. In particular, the $\mathcal{S}^*(\alpha)$-radius of a normalized function $f$ is to find
	$$\sup\{ r\in \mathbb{R}^+ : \RE\left(\frac{zf'(z)}{f(z)} \right) > \alpha, z\in \mathbb{D}_{r}  \}$$
	and similarly, we can define $\mathcal{C}(\alpha)$-radius. In this direction, for $\mathcal{S}^*(\alpha)$-radius and more, we refer to Bessel functions~\cite{abo-2018,Baricz,Ba3} (see Watson's treatise~\cite{watson-1944} for more on Bessel function), Struve functions~\cite{abo-2018,bdoy-2016}, Wright functions~\cite{btk-2018}, Lommel functions~\cite{abo-2018,bdoy-2016} and Ramanujan type entire functions~\cite{ErhanDenij2020}. For their generalization to Ma-Minda classes~\cite{minda94} of starlike and convex functions, we refer to see~\cite{g-specialIJST,SG-2020}.

	However, with the best of our knowledge, $\mathcal{S}_{p}^{\gamma}(\alpha)$-radius and $\mathcal{CS}_{p}^{\gamma}(\alpha)$-radius for special functions are not handled till date. Therefore, we define $\mathcal{S}_{p}^{\gamma}(\alpha)$-radius and its convex analog here below: 
	\begin{definition}
	   Let  $f$ in $\mathcal{A}$ be a special function. Then the radius of $\gamma$-Spirallike of order $\alpha$ is to find
		\begin{equation*}
		r^*_{sp}(\alpha,\gamma;f) =\sup\left\{ r\in \mathbb{R}^+ : \RE\left(e^{-i\gamma} \frac{zf'(z)}{f(z)} \right) > \alpha \cos{\gamma}, z\in \mathbb{D}_{r} \right \}
		\end{equation*}
		and the radius of convex $\gamma$-Spirallike of order $\alpha$ is given by
		\begin{equation*}
		r_{sp}^{c}(\alpha,\gamma;f) = \sup\left\{ r\in \mathbb{R}^+ : \RE\left(e^{-i\gamma} \left(1+\frac{zf''(z)}{f'(z)}\right) \right) > \alpha \cos{\gamma}, z\in \mathbb{D}_{r} \right \}.
		\end{equation*}
	\end{definition}

	In the present investigation, we consider our special normalized functions to be the derivatives of Bessel functions. For this to proceed further, recall that the Bessel function of the first kind of order $\nu $ is defined by \cite[p.217]{Olver} 
\begin{equation}
J_{\nu }(z)=\sum_{n=0}^{\infty }\frac{\left( -1\right) ^{n}}{n!\Gamma (n+\nu+1)}\left( \frac{z}{2}\right) ^{2n+\nu }\text{ \ \ }
\left( z\in 
\mathbb{C}
\right) .  \label{J1}
\end{equation}%
We know that it has all its zeros real for $\nu >-1$. Here now we consider mainly the general function
\begin{equation*}
N_{\nu }(z)=az^{2}J_{\nu }^{\prime \prime }(z)+bzJ_{\nu }^{\prime}(z)+cJ_{\nu }(z)
\end{equation*}%
studied by Mercer \cite{Mercer}. Here, as in \cite{Mercer}, $q=b-a$ and 
\begin{equation*}
\left(c=0 \text{ and } q\neq0 \right) \text{ or } \left(c>0 \text{ and } q>0 \right) .
\end{equation*}
From (\ref{J1}), we have the power series representation%
\begin{equation}
N_{\nu }(z)=\sum_{n=0}^{\infty }\frac{Q(2n+\nu )\left( -1\right) ^{n}}{n!\Gamma (n+\nu +1)}\left( \frac{z}{2}\right) ^{2n+\nu }\text{ \ \ }
\left(
z\in 
\mathbb{C}
\right)  \label{J2}
\end{equation}%
where $Q(\nu)=a\nu(\nu-1)+b\nu+c$ $\left( a,b,c\in 
\mathbb{R}
\right) .$ There are three important works on the function $N_{\nu }.$ Firstly, Mercer's paper \cite{Mercer} in which it has been proved that the $k^{th}$ positive zero of $N_{\nu }$ increases with $\nu $ in $\nu >0$. Secondly, Ismail and Muldoon \cite{IS} showed that under the conditions $a,b,c\in 
\mathbb{R}
$ such that $c=0$ and $b\neq a$ or $c>0$ and $b>a$;
\begin{itemize}
	\item[(i)] For $\nu > 0$, the zeros of $N_{\nu }(z)$ are either real or
	purely imaginary.
	
	\item[(ii)] For $\nu \geq \max \{0,\nu _{0}\}$, where $\nu _{0}$ is the
	largest real root of the quadratic $Q(\nu )=a\nu (\nu -1)+b\nu +c,$ the the
	zeros of $N_{\nu }(z)$ are real.
	
	\item[(iii)] If $\nu >0$, $\left( a\nu ^{2}+(b-a)\nu +c\right) \diagup
	(b-a)>0$ and $a\diagup (b-a)<0$, the zeros of $N_{\nu }(z)$ are all real
	except for a single pair which are conjugate purely imaginary.
\end{itemize}
Baricz et al.~\cite{Erhandeniz} obtained sufficient and necessary conditions for the starlikeness of a normalized form of $N_{\nu }$ by using results of Mercer~\cite{Mercer}, Ismail and Muldoon~\cite{IS} and Shah and Trimble~\cite{Shah}. 

Note that $N_{\nu }$ do not belong to $\mathcal{A}$. To prove our main results, we consider the following normalizations of the function $N_{\nu }$ given by:
\begin{eqnarray}
f_{\nu }(z) &=&\left[ \frac{2^{\nu }\Gamma (\nu +1)}{Q(\nu)}N_{\nu }(z)%
\right] ^{\frac{1}{\nu }},  \label{eq13} \\
g_{\nu }(z) &=&\frac{2^{\nu }\Gamma (\nu +1)z^{1-\nu }}{Q(\nu)}N_{\nu }(z), 
\label{eq14} \\
h_{\nu }(z) &=&\frac{2^{\nu }\Gamma (\nu +1)z^{1-\frac{\nu }{2}}}{Q(\nu)}%
N_{\nu }(\sqrt{z}).  \label{eq15}
\end{eqnarray}
In the rest of this paper, for the quadratic $Q(\nu )=a\nu (\nu -1)+b\nu +c$, we will always assume that $a,b,c\in \mathbb{R} $ $\left(c=0 \text{ \  and \ } a\neq b    \right)$ or $\left(c> 0 \text{ \  and \ } a< b    \right).$ Moreover, $\nu_0$ is the largest real root of the quadratic $Q(\nu )$ defined according to the above conditions.

One can see that the functions $N_\nu(z)$ and $N^\prime_\nu(z)$
are entire functions of order zero, then they have infinitely many zeros. According to Hadamard factorization theorem \cite{Le}
\begin{equation}\label{NVZCrp}
N_{\nu }(z)=\frac{Q(\nu )z^{\nu }}{2^{\nu }\Gamma (\nu +1)}\prod\limits_{n\geq 1}\left(
1-\frac{z^{2}}{ \lambda_{\nu ,n} ^{2}}\right)
\end{equation}
and
\begin{equation}\label{NVZTCrp}
N^\prime_{\nu }(z)=\frac{Q(\nu ) z^{\nu-1 } \nu}{2^{\nu }\Gamma (\nu +1)}\prod\limits_{n\geq
	1}\left( 1-\frac{z^{2}}{ \lambda_{\nu ,n} ^{\prime2}}\right)
\end{equation}
where $\lambda_{\nu ,n}$ and $\lambda_{\nu ,n} ^{\prime}$ denote the $ n $th positive zero of $N_\nu(z)$ and $N^\prime_\nu(z),$ respectively. Recently, for studies on the Geometric properties of Bessel functions, see \cite{Baricz,Ba1,Ba3,Erhandeniz,Ba21,Br,Ca,De,Deniztrev,Kr,Rav,Sza}.

Recall that Ma and Minda \cite{minda94} subclasses of starlike and convex functions, respectively, are given by: 
$$\mathcal{S}^*\left(\varphi\right)=\left\{f\in \mathcal{A}:\frac{zf^\prime(z)}{f(z)}\prec  \varphi\left(z\right) \right\}  $$
and
$$\mathcal{C}\left(\varphi\right)=\left\{f\in \mathcal{A}:1+\frac{zf^{\prime\prime}(z)}{f^\prime(z)}\prec  \varphi\left(z\right) \right\},  $$
where $ \varphi $, a Ma-Minda function, is analytic and univalent with $\operatorname{Re}\varphi\left(z\right)>0,~\varphi^\prime\left(0\right)>0$ and $ \varphi\left(\mathbb{D}\right) $ is starlike with respect to $\varphi\left(0\right)=1$ and symmetric about real axis. Note that $\varphi\in \mathcal{P},$ the class of normalized Carath\'{e}odory functions. 	
We also recall that
	\begin{definition}
		Let $f\in \mathcal{A}$ be a special function. Then $\mathcal{S}^*\left(\varphi\right)-$radius and $\mathcal{C}\left(\varphi\right)-$radius of $ f $ are defined as follows:
		$$r^*_\varphi\left(f\right)=\sup\left\{r\in \mathbb{R}^+ :\frac{zf^\prime(z)}{f(z)}\in  \varphi\left(\mathbb{D}\right),~z\in \mathbb{D}_{r}  \right\}    $$
		and
		$$r^c_\varphi\left(f\right)=\sup\left\{r\in \mathbb{R}^+ :1+\frac{zf^{\prime\prime}(z)}{f^\prime(z)}\in  \varphi\left(\mathbb{D}\right),~z\in \mathbb{D}_{r}  \right\}, $$
		respectively.
	\end{definition}

	In this paper, we deal with the radius of $\gamma$-Spirallike of order $\alpha$ and the radius of convex $\gamma$-Spirallike of order $\alpha$ for the functions $f_{\nu }(z),\;g_{\nu }(z)$ and $h_{\nu }(z)$ defined by \eqref{eq13}, \eqref{eq14} and \eqref{eq15} in the case when $\nu \geq \max \{0,\nu _{0}\}$. Also, we determine the $\mathcal{S}^*\left(\varphi\right)-$radii and $\mathcal{C}\left(\varphi\right)-$radii of these functions. The key tools in their proofs are some new Mittag-Leffler expansions for quotients of the function $N_{\nu }$, special properties of the zeros of the function $N_{\nu }$ and their derivatives.

	\section{Zeros of hyperbolic polynomials and the Laguerre--P\'{o}lya class of entire functions}
	
	In this section, we recall some necessary information about polynomials and entire functions with real zeros. An algebraic polynomial is called hyperbolic if all its zeros are real. We formulate the following specific statement that we shall need, see \cite{Deniztrev} for more details.
	
	By definition, a real entire function $\psi $ belongs to the Laguerre--P\'{o}lya class $\mathcal{LP}$ if it can be represented in the form
	\begin{equation*}
	\psi (x)=cx^{m}e^{-ax^{2}+\beta x}\prod\limits_{k\geq 1}\left( 1+\frac{x}{x_{k}}\right) e^{-\frac{x}{x_{k}}},
	\end{equation*}%
	with $c,\beta ,x_{k}\in 
	\mathbb{R}
	,$ $a\geq 0,$ $m\in 
	\mathbb{N}
	\cup \{0\}$ and $\sum x_{k}^{-2}<\infty .$ Similarly, $\phi $ is said to be of type $\mathcal{I}$ in the Laguerre-P\'{o}lya class, written $\varphi \in 
	\mathcal{LPI}$, if $\phi (x)$ or $\phi (-x)$ can be represented as%
	\begin{equation*}
	\phi (x)=cx^{m}e^{\sigma x}\prod\limits_{k\geq 1}\left( 1+\frac{x}{x_{k}}\right) ,
	\end{equation*}%
	with $c\in 
	\mathbb{R}
	,$ $\sigma \geq 0,$ $m\in 
	\mathbb{N}
	\cup \{0\},$ $x_{k}>0$ and $\sum x_{k}^{-1}<\infty .$ The class $\mathcal{LP}$ is the complement of the space of hyperbolic polynomials in the topology induced by the uniform convergence on the compact sets of the complex plane while $\mathcal{LPI}$ is the complement of the hyperbolic polynomials whose zeros possess a preassigned constant sign. Given an entire function $\varphi $ with the Maclaurin expansion%
	\begin{equation*}
	\varphi (x)=\sum_{k\geq 0}\mu _{k}\frac{x^{k}}{k!},
	\end{equation*}%
	its Jensen polynomials are defined by 
	\begin{equation*}
	P_{m}(\varphi ;x)=P_{m}(x)=\sum_{k=0}^{m}\left( 
	\begin{array}{c}
	m \\ 
	k%
	\end{array}%
	\right) \mu _{k}x^{k}.
	\end{equation*}%
	The next result of Jensen \cite{Je} is a well-known characterization of functions belonging to $\mathcal{LP}$.
	
	\begin{lemma}
		\label{Le2} The function $\varphi $ belongs to $\mathcal{LP}$ $(\mathcal{LPI}$, respectively$)$ if and only if all the polynomials $P_{m}(\varphi ;x)$, $m=1,2,...,$ are hyperbolic (hyperbolic with zeros of equal sign). Moreover, the sequence $P_{m}(\varphi ;z\diagup n)$ converges locally uniformly to $\varphi (z)$.
	\end{lemma}
	
	The following result is a key tool in the proof of main results.
	
	\begin{lemma}\cite{Kazimoglu}\label{Le4}
		If $\nu \geq \max\{0,\nu_0\}$ then the functions $z\longmapsto\Psi_\nu(z)=\frac{2^{\nu }\Gamma (\nu +1)}{Q(\nu )z^{\nu }}N_\nu(z)$ has infinitely many zeros and all of them are positive. Denoting by $\lambda_{\nu,n}$ the $n$th positive zero of $\Psi_\nu(z),$ under the same conditions the Weierstrassian decomposition 
		\begin{equation*}
		\Psi_\nu(z)=\prod\limits_{n\geq 1}\left( 1-\frac{z^{2}}{ \lambda_{\nu ,n}^{2}}\right)
		\end{equation*}
		is valid, and this product is uniformly convergent on compact subsets of the complex plane. Moreover, if we denote by $\lambda_{\nu ,n}^\prime$ the $n$th positive zero of $\Phi_\nu^\prime(z),$ where $\Phi_\nu(z)=z^\nu \Psi_\nu(z), $ then the positive zeros of $\Psi_\nu(z)$ are interlaced with those of $\Phi_\nu^\prime(z).$ In the other words, the zeros satisfy the chain of inequalities 
		\begin{equation*}
		\lambda_{\nu ,1}^\prime<\lambda_{\nu ,1}<\lambda_{\nu,2}^{\prime}<\lambda_{\nu ,2}<\lambda_{\nu ,3}^\prime<\lambda_{\nu,3}<\cdots .
		\end{equation*}
	\end{lemma}

	\section{Main \ Results}
	
	\subsection{Radii of $\gamma-$Spirallike and convex $\gamma-$Spirallike of The Functions $f_{ \nu },$ $g_{\nu }$ and $h_{\nu }$} 
	The first principal result we established concerns the radii of $\gamma-$Spirallike of order $\alpha$ and reads as follows.
	
	\begin{theorem}\label{T1} 
		Let $\alpha \in\lbrack 0,1)$ and $\gamma \in \left(-\frac{\pi}{2},\frac{\pi}{2}\right).$ The following statements hold:
		
		\begin{enumerate}
			\item[a)] If $\nu \geq \max\{0,\nu_{0} \}$, $\nu\neq0$  then the radius of $\gamma-$Spirallike of order $\alpha$ of the function $f_\nu$ is the smallest positive root of the equation 
			\begin{equation*}
			\frac{1}{\nu}\frac{ar^3J_\nu^{\prime\prime\prime}(r)+\left(2a+b\right)r^2J_\nu^{\prime\prime}(r)+\left(b+c\right)rJ_\nu^{\prime}(r)}{%
				ar^2J_\nu^{\prime\prime}(r)+brJ_\nu^{\prime}(r)+cJ_\nu(r)}=1-\left(1-\alpha\right)\cos \gamma.
			\end{equation*}
			
			\item[b)] If $\nu \geq \max\{0,\nu_{0} \},$ then the radius of $\gamma-$Spirallike of order $\alpha$ of the function $g_\nu$ is the smallest positive root of the equation 
			\begin{equation*}
			\frac{ar^3J_\nu^{\prime\prime\prime}(r)+\left(2a+b%
				\right)r^2J_\nu^{\prime\prime}(r)+\left(b+c\right)rJ_\nu^{\prime}(r)}{%
				ar^2J_\nu^{\prime\prime}(r)+brJ_\nu^{\prime}(r)+cJ_\nu(r)}=\nu-\left(1-\alpha\right)\cos \gamma.
			\end{equation*}
			
			\item[c)] If $\nu \geq \max\{0,\nu_{0} \},$ then the radius of $\gamma-$Spirallike of order $\alpha$ of the function $%
			h_\nu$ is the smallest positive root of the equation 
			\begin{equation*}
			\frac{ar\sqrt{r} J_\nu^{\prime\prime\prime}(\sqrt{r}%
				)+\left(2a+b\right)rJ_\nu^{\prime\prime}(\sqrt{r})+\left(b+c\right)\sqrt{r}%
				J_\nu^{\prime}(\sqrt{r})}{arJ_\nu^{\prime\prime}(\sqrt{r})+b\sqrt{r}%
				J_\nu^{\prime}(\sqrt{r})+cJ_\nu(\sqrt{r})}=\nu-2\left(1-\alpha\right)\cos \gamma.
			\end{equation*}
		\end{enumerate}
	\end{theorem}
	
	\begin{proof}
		Firstly, we prove part \textbf{a} for $\nu \geq \max \{0,\nu_0\},~\nu\neq0 $
		and \textbf{b} and \textbf{c} for $\nu \geq \max \{0,\nu_0\}.$ We need to
		show that the following inequalities for $\alpha \in
		\lbrack 0,1)$ and $\gamma \in \left(-\frac{\pi}{2},\frac{\pi}{2}\right),$
		\begin{equation} \small
		\operatorname{Re}\left( e^{-i\gamma}\frac{zf_{\nu }^{\prime }(z)}{f_{\nu }(z)}\right) >\alpha\cos\gamma,
		~~~ \operatorname{Re}\left( e^{-i\gamma}\frac{zg_{\nu }^{\prime }(z)}{g_{\nu }(z)}\right) >\alpha\cos\gamma
		\text{ and } 
		\operatorname{Re}\left(e^{-i\gamma} \frac{zh_{\nu }^{\prime }(z)}{h_{\nu }(z)}\right) >\alpha\cos\gamma  \label{eq21}
		\end{equation}\normalsize
		are valid for $z\in {\mathbb{D}}_{r_{sp }^{*}(\alpha,\gamma;f_{\nu })},~z\in \mathbb{D}_{r_{sp }^{*}(\alpha,\gamma;g_{\nu })}$ and $z\in \mathbb{D}_{r_{sp }^{*}(\alpha,\gamma;h_{\nu })}$ respectively, and each of the above inequalities does not hold
		in larger disks.\\
		\textbf{a):} When we write the equation (\ref{eq13}) in definition of the function $f_{\nu }(z),$ we get by using logarithmic derivation 
		\begin{equation}
		\frac{zf_{\nu }^{\prime }(z)}{f_{\nu }(z)} =\frac{1}{\nu }\frac{ zN_{\nu
			}^\prime(z)}{N_{\nu }(z)}=1-\frac{1}{\nu}\sum\limits_{n\geq 1}\frac{2z^{2}}{
			\lambda_{\nu ,n} ^{2}-z^{2} },\text{ \ \ }
		\left(\nu \geq \max\{0,\nu_0\},~\nu\neq0\right),\label{eq22} 
		\end{equation}
		It is known \cite{De} that if $z\in \mathbb{C}$ and $\lambda \in \mathbb{R}$ are such that $ \left\vert z\right\vert\leq r<\lambda ,$ then
		\begin{equation}
		\operatorname{Re}\left( \frac{z}{\lambda -z}\right) \leq \left\vert \frac{z}{\lambda -z}\right\vert  \leq \frac{|z|}{\lambda -\left\vert z\right\vert }.  \label{eq23}
		\end{equation}%
		Then the inequality
		\begin{equation*}
		\operatorname{Re}\left( \frac{z^{2}}{ \lambda_{\nu,n}^{2}-z^{2}}\right) \leq \left\vert \frac{z^{2}}{ \lambda_{\nu,n}^{2}-z^{2}}\right\vert \leq \frac{\left\vert z\right\vert ^{2}}{ \lambda_{\nu ,n}^ {2}-\left\vert z\right\vert ^{2}}
		\end{equation*}%
		holds for every $\left|z\right|<\lambda_{\nu ,1}. $ Therefore, from (\ref{eq22}) and (\ref{eq23}), we have 
		\begin{eqnarray*}
			\operatorname{Re}\left(  e^{-i\gamma}\frac{zf_{\nu }^{\prime }(z)}{f_{\nu }(z)}\right) &=& \operatorname{Re}\left( e^{-i\gamma}\right)-\frac{1}{\nu}\operatorname{Re}\left(e^{-i\gamma}\sum\limits_{n\geq 1} \frac{2z^{2}}{ \lambda_{\nu
					,n}^{2}-z^{2}}\right)  \\
			&\geq& \cos\gamma- \frac{1}{\nu}\left|e^{-i\gamma}\sum\limits_{n\geq 1} \frac{2z^{2}}{ \lambda_{\nu
					,n}^{2}-z^{2}}\right|\geq  \cos\gamma-\frac{1}{\nu }\sum\limits_{n\geq 1}\frac{%
				2\left\vert z\right\vert ^{2}}{ \lambda_{\nu ,n}^{2}-\left\vert z\right\vert
				^{2}}    \\
			&=&\frac{\left\vert z\right\vert f_{\nu }^{\prime }(\left\vert
				z\right\vert )}{f_{\nu }(\left\vert z\right\vert )}+\cos\gamma-1,
		\end{eqnarray*}%
		with equality when $z=\left\vert z\right\vert =r.$
		Thus, for $r\in \left(0,\lambda_{\nu ,1} \right)$ it follows that 
		\begin{equation*}
		\inf_{z\in \mathbb{D}_r } \left\lbrace \operatorname{Re} \left( e^{-i\gamma}\frac{zf_{\nu }^{\prime }(z)}{f_{\nu }(z)}-\alpha\cos\gamma \right) \right\rbrace= \frac{\left\vert z\right\vert f_{\nu }^{\prime }(\left\vert
			z\right\vert )}{f_{\nu }(\left\vert z\right\vert )}+\left(1-\alpha\right)\cos\gamma-1.
		\end{equation*}
		Now, the mapping $\Theta_\nu:\left(0,\lambda_{\nu ,1} \right)\longrightarrow \mathbb{R}$ defined by 
		\begin{equation*}
		\Theta_\nu(r)= \frac{rf_{\nu }^{\prime }(r)}{f_{\nu }(r)}+\left(1-\alpha\right)\cos\gamma-1=\left(1-\alpha\right)\cos\gamma-\frac{1}{\nu}%
		\sum\limits_{n\geq 1}\left( \frac{2r^{2}}{ \lambda_{\nu ,n}^{2}-r^{2}}%
		\right).
		\end{equation*}
		is strictly decreasing since
		\begin{equation*}
		\Theta_\nu^\prime(r)= -\frac{1}{\nu}\sum\limits_{n\geq 1}\left( \frac{%
			4r\lambda_{\nu,n}^{2}}{ \left(\lambda_{\nu ,n}^{2}-r^{2}\right)^2}\right)<0
		\end{equation*}
		for all $\nu\geq
		\max\{0,\nu_0 \}, \nu \neq 0.$ On the other hand, since $$\lim_{r\searrow0}\Theta_\nu(r)=\left(1-\alpha\right)\cos\gamma>0  \text{ \  and \ } \lim_{r\nearrow \lambda_{\nu ,1}}\Theta_\nu(r)=-\infty,$$
		in view of the minimum principle for harmonic functions imply that the corresponding inequality in (\ref{eq21})
		for $\nu\geq \max\{0,\nu_0 \}, \nu \neq 0$ holds if and only if $z\in \mathbb{D}_{r_f},$ where $r_f$ is the smallest positive root of equation 
		\begin{equation*}
		\frac{rf_{\nu }^{\prime }(r)}{f_{\nu }(r)}=1-\left(1-\alpha\right)\cos\gamma
		\end{equation*}
		which is equivalent to 
		\begin{equation*}
		\frac{rN_\nu^\prime(r)}{\nu N_\nu(r)}=1-\left(1-\alpha\right)\cos\gamma, 
		\end{equation*}
		situated in $\left(0,\lambda_{\nu ,1}\right).$\\
		\textbf{b):} Since,
		\begin{equation} \label{eq24}
		\frac{zg_{\nu }^{\prime }(z)}{g_{\nu }(z)} =(1-\nu) +\frac{zN_{\nu}^{\prime}(z)}{N_{\nu }(z)}=1- \sum\limits_{n\geq 1}\frac{2z^{2}}{\lambda_{\nu ,n}^2-z^{2}},
		\text{ \ \ } \left(\nu \geq \max\{0,\nu_0\}\right).
		\end{equation}
		By using the inequality (\ref{eq23}), for all $z \in \mathbb{D}_{\lambda_{\nu ,1} }$ we get the inequality  
		\begin{eqnarray*}
			\operatorname{Re}\left(  e^{-i\gamma}\frac{zg_{\nu }^{\prime }(z)}{g_{\nu }(z)}\right) 
			&=& \operatorname{Re}\left( e^{-i\gamma}\right)-\operatorname{Re}\left(e^{-i\gamma}\sum\limits_{n\geq 1} \frac{2z^{2}}{ \lambda_{\nu,n}^{2}-z^{2}}\right)  \\
			&\geq& \cos\gamma- \left|e^{-i\gamma}\sum\limits_{n\geq 1} \frac{2z^{2}}{ \lambda_{\nu,n}^{2}-z^{2}}\right|\geq  \cos\gamma-\sum\limits_{n\geq 1}\frac{2\left\vert z\right\vert ^{2}}{ \lambda_{\nu ,n}^{2}-\left\vert z\right\vert^{2}}    \\
			&=&\frac{r g_{\nu }^{\prime }(r)}{g_{\nu }(r)}+\cos\gamma-1,
		\end{eqnarray*}%
		where $\left\vert z\right\vert =r.$
		Thus, for $r\in \left(0,\lambda_{\nu ,1} \right)$ we obtain
		\begin{equation*}
		\inf_{z\in \mathbb{D}_r } \left\lbrace \operatorname{Re} \left( e^{-i\gamma}\frac{zg_{\nu }^{\prime }(z)}{g_{\nu }(z)}-\alpha\cos\gamma \right) \right\rbrace= \frac{r g_{\nu }^{\prime }(r)}{g_{\nu }(r)}+\left(1-\alpha\right)\cos\gamma-1.
		\end{equation*}
		The function $\Phi_\nu:\left(0,\lambda_{\nu ,1} \right)\longrightarrow \mathbb{R}$ defined by 
		\begin{equation*}
		\Phi_\nu(r)= \frac{rg_{\nu }^{\prime }(r)}{g_{\nu }(r)}+\left(1-\alpha\right)\cos\gamma-1=\left(1-\alpha\right)\cos\gamma-\sum\limits_{n\geq 1}\left( \frac{2r^{2}}{ \lambda_{\nu,n}^{2}-r^{2}}%
		\right).
		\end{equation*}
		is strictly decreasing and $\lim_{r\searrow0}\Phi_\nu(r)=\left(1-\alpha\right)\cos\gamma>0$ and $\lim_{r\nearrow \lambda_{\nu ,1}}\Phi_\nu(r)=-\infty.$\\ Consequently, in view of the minimum principle for harmonic functions for $z\in \mathbb{D}_{r_g},$ we have that
		$$ \operatorname{Re} \left( e^{-i\gamma}\frac{zg_{\nu }^{\prime }(z)}{g_{\nu }(z)}\right)>\alpha\cos\gamma  $$
		if and only if $r_g$ is the smallest positive root of equation 
		\begin{equation*}
		\frac{rg_{\nu }^{\prime }(r)}{g_{\nu }(r)}=1-\left(1-\alpha\right)\cos\gamma \text{ \  or \ } \frac{rN_\nu^\prime(r)}{ N_\nu(r)}=\nu-\left(1-\alpha\right)\cos\gamma
		\end{equation*}
		situated in $\left(0,\lambda_{\nu ,1}\right).$\\
		\textbf{c):} Observe that
		\begin{equation}
		\frac{zh_{\nu }^{\prime }(z)}{h_{\nu }(z)}=(1-\frac{\nu}{2}) +\frac{1}{2} 
		\frac{\sqrt{z} N_{\nu }^{\prime}(\sqrt{z})}{N_{\nu }(\sqrt{z})}=1-
		\sum\limits_{n\geq 1}\frac{z}{ \lambda_{\nu ,n}^2-z},\text{ \ \ }\left(\nu \geq
		\max\{0,\nu_0 \}\right),\label{eq25} 
		\end{equation}
		By using the inequality (\ref{eq22}), for all $z\in \mathbb{D}_{\lambda_{\nu ,1}^2}$ we obtain
		\begin{eqnarray*}
			\operatorname{Re}\left(  e^{-i\gamma}\frac{zh_{\nu }^{\prime }(z)}{h_{\nu }(z)}\right) &=& \operatorname{Re}\left(e^{-i\gamma}\right)-\operatorname{Re}\left(e^{-i\gamma}\sum\limits_{n\geq 1} \frac{z}{ \lambda_{\nu,n}^{2}-z}\right)  \\
			&\geq& \cos\gamma-\left|e^{-i\gamma}\sum\limits_{n\geq 1} \frac{z}{ \lambda_{\nu,n}^{2}-z}\right|\geq  \cos\gamma-\sum\limits_{n\geq 1}\frac{%
				2\left\vert z\right\vert }{ \lambda_{\nu ,n}^{2}-\left\vert z\right\vert}    \\
			&=&\frac{r h_{\nu }^{\prime }(r)}{h_{\nu }(r )}+\cos\gamma-1,
		\end{eqnarray*}%
		Equality holds for $z=\left\vert z\right\vert =r.$
		As a result, for $r\in \left(0,\lambda_{\nu ,1} \right)$ it follows that 
		\begin{equation*}
		\inf_{z\in \mathbb{D}_r } \left\lbrace \operatorname{Re} \left( e^{-i\gamma}\frac{zh_{\nu }^{\prime }(z)}{h_{\nu }(z)}-\alpha\cos\gamma \right) \right\rbrace= \frac{r h_{\nu }^{\prime }(r )}{h_{\nu }(r )}+\left(1-\alpha\right)\cos\gamma-1.
		\end{equation*}
		The fucntion $\Psi_\nu:\left(0,\lambda_{\nu ,1}^2 \right)\longrightarrow \mathbb{R}$ defined by 
		\begin{equation*}
		\Psi_\nu(r)= \frac{rh_{\nu }^{\prime }(r)}{h_{\nu }(r)}+\left(1-\alpha\right)\cos\gamma-1=\left(1-\alpha\right)\cos\gamma-\sum\limits_{n\geq 1}\left( \frac{2r}{ \lambda_{\nu ,n}^{2}-r}%
		\right).
		\end{equation*}
		is strictly decreasing, since
		\begin{equation*}
		\Psi_\nu^\prime(r)= -\sum\limits_{n\geq 1}\left( \frac{\lambda_{\nu,n}^{2}}{ \left(\lambda_{\nu ,n}^{2}-r\right)^2}\right)<0
		\end{equation*}
		for all $\nu\geq\max\{0,\nu_0 \}.$ 
		On the other hand, since $$\lim_{r\searrow0}\Psi_\nu(r)=\left(1-\alpha\right)\cos\gamma>0  \text{ \  and \ } \lim_{r\nearrow \lambda_{\nu ,1}^2}\Psi_\nu(r)=-\infty$$
		in view of the minimum principle for
		harmonic functions imply that the corresponding inequality in (\ref{eq21}) holds if and only if $z\in \mathbb{D}_{r_h},$ where $r_h$ is the smallest positive root of equation 
		\begin{equation*}
		\frac{rh_{\nu }^{\prime }(r)}{h_{\nu }(r)}=1-\left(1-\alpha\right)\cos\gamma
		\end{equation*}
		which is equivalent to 
		\begin{equation*}
		\frac{\sqrt{r}N_\nu^\prime(\sqrt{r})}{N_\nu(\sqrt{r})}=1-\left(1-\alpha\right)\cos\gamma 
		\end{equation*}
		and situated in $\left(0,\lambda_{\nu ,1}^2\right).$ This completes the proof of part \textbf{a} when $\nu\geq \max\{0,\nu_0 \},
		\nu \neq 0$, and parts \textbf{b} and \textbf{c} when $\nu\geq \max\{0,\nu_0\}.$ \qed
	\end{proof} 

\begin{remark}	
	In Theorem~\ref{T1}, the choice of $\gamma=0$ yields the result for the class of starlike functions of order $\alpha$.
\end{remark}

\begin{example}
	Indeed, we can write the function $N_{\nu}(z)$ in terms of
	elementary trigonometric functions for the value $\nu=1/2$ as follows:
	\begin{equation*}
	N_{1/2}(z)=\frac{4\left(b-a\right)z\cos z+\left[a\left( 3-4z^{2}\right)-2b+4c\right] \sin z}{2\sqrt{2\pi }\sqrt{z%
	}}.
	\end{equation*}
	Thus, we have
	\begin{equation*}
	f_{1/2}(z)=\frac{\left[ 4\left(a-b\right)z\cos z+\left( 4az^{2}-3a+2b-4c\right) \sin z\right] ^{2}}{\left(a-2b-4c\right)^2z},
	\end{equation*}%
	\begin{equation*}
	g_{1/2}(z)=\frac{4\left(a-b\right)z\cos z+\left( 4az^{2}-3a+2b-4c\right) \sin z}{a-2b-4c}
	\end{equation*}
	and 
	\begin{equation*}
	h_{1/2}(z)=\frac{4\left(a-b\right)z\cos \sqrt{z}+\left( 4az-3a+2b-4c\right)\sqrt{z} \sin \sqrt{z}}{a-2b-4c}.
	\end{equation*}
\end{example}		
	
	Taking $\nu=1/2,$ $a=1,$ $b=2,$ $c=0,$ $\alpha=1/2$ and $\gamma=\pi/3$ in Theorem \ref{T1}, we have following results.
	\begin{corollary}\label{corlySprlke} The following statements are true.
		\begin{itemize}
			\item[\textbf{a)}] The radius of $\frac{\pi}{3}-$Spirallike of order $1/2$ of the function $f_{1/2}$ is the smallest positive root of the equation
			$$ \frac{4x\left(8x+1\right)\cos{x}+\left(68x^2-7\right)\sin{x}}{4\left(-4x\cos{x}+\left(4x^2+1\right)\sin{x}\right)}=0$$
			
			\item[\textbf{b)}] The radius of $\frac{\pi}{3}-$Spirallike of order $1/2$ of the function $g_{1/2}$ is the smallest positive root of the equation
			$$ \frac{16x^3\cos{x}+3\left(12x^2-1\right)\sin{x}}{4\left(-4x\cos{x}+\left(4x^2+1\right)\sin{x}\right)}=0$$
			
			\item[\textbf{c)}] The radius of $\frac{\pi}{3}-$Spirallike of order $1/2$ of the function $h_{1/2}$ is the smallest positive root of the equation
			$$ \frac{2\sqrt{x}\left(4x-1\right)\cos{\sqrt{x}}+\left(20x-1\right)\sin{\sqrt{x}}}{2\left(-4\sqrt{x}\cos{\sqrt{x}}+\left(4x+1\right)\sin{\sqrt{x}}\right)}=0.$$
			
		\end{itemize}
	\end{corollary}

	\begin{table} 
		\caption{ $R\left(--\right):$ Radii of $\frac{\pi}{3}-$Spirallike of order $1/2$ for $f_{1/2 },~g_{1/2}$ and $h_{1/2 }$.}
		\begin{tabular}{|c|c|c|c|c|c|c|c|c|c|}
			\hline
			& \multicolumn{3}{c|}{$b=3$ and $c=0$} & \multicolumn{3}{c|}{$a=1$ and $c=0$}
			& \multicolumn{3}{c|}{$a=1$ and $b=2$} \\ \cline{2-10}
			& $a=2$ & $a=3$ & $a=4$ & $b=2$ & $b=3$ & $b=4$ & $c=2$ & $c=3$ & $c=4$ \\ 
			\hline
			\multicolumn{1}{|c|}{  $R\left(f_{\frac{1}{2}}\right)  $} & $0.1539$ & $0.1190$
			& $0.0886$ & $0.1746$ & $0.1990$ & $0.2130$ & $0.3038$ & $0.3401$ & $0.3680$
			\\ \hline
			\multicolumn{1}{|c|}{$R\left(g_{\frac{1}{2}}\right) $} & $0.2121$ & $0.1639$
			& $0.1220$ & $0.2409$ & $0.2747$ & $0.2942$ & $0.4217$ & $0.4730$ & $0.5126$
			\\ \hline
			\multicolumn{1}{|c|}{$R\left(g_{\frac{1}{2}}\right) $} & $0.0817$ & $0.0486$
			& $0.0268$ & $0.1056$ & $0.1377$ & $0.1582$ & $0.3309$ & $0.4193$ & $0.4953$
			\\ \hline
		\end{tabular}
		 \label{Table1.}
	\end{table}

\begin{figure}[h]
	\begin{tabular}{c}
		\includegraphics[scale=0.9]{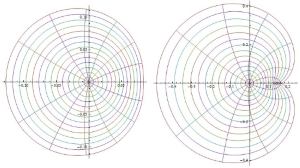}
	\end{tabular}
	\caption{Images of function $h_{1/2}(z)$ for $r=0.1056$ and $r=0.3$, respectively as per  Corollary~\ref{corlySprlke}.}\label{Figure1.}
\end{figure}

The figure~\ref{Figure1.} was created by taking the special values in the Corollary~\ref{corlySprlke}.

\begin{remark}
	When $\nu=0.5,$ considering the special values of $a,~b$ and $c$ real numbers, radii of $\gamma-$Spirallike of order $\alpha$ of the function $f_\nu$ is seen from the table above. In addition, according to the increasing values of $a,~b$ and $c$ reel numbers, it is clear that radii of $\gamma-$Spirallike of order $\alpha$ of the function $f_\nu$ is monotone.
\end{remark}

	In the following, we deal with convex analogue of the class of $\gamma$-spiralllike functions of order $\alpha$.	
	\begin{theorem}\label{T2} 
		Let $\alpha \in\lbrack 0,1)$ and $\gamma \in \left(-\frac{\pi}{2},\frac{\pi}{2}\right).$ The following statements hold:	
		\begin{enumerate}
			\item[a)] If $\nu \geq \max\{0,\nu_{0} \}$, $\nu\neq0$  then the radius of the convex $\gamma-$Spirallike of order $\alpha$ of the function $f_\nu$ is the smallest positive root of the equation 
			\begin{equation*}
			\frac{rN_{\nu}^{\prime \prime }(r)}{N_{\nu }^{\prime }(r)}+\left(\frac{%
				1}{\nu}-1\right)\frac{rN_{\nu}^{ \prime }(r)}{N_{\nu }(r)}=\left(\alpha-1\right)\cos\gamma.
			\end{equation*}
			
			\item[b)] If $\nu \geq \max\{0,\nu_{0} \},$ then the radius of the convex $\gamma-$Spirallike of order $\alpha$ of the function $g_\nu$ is the smallest positive root of the equation 
			\begin{equation*}
			\frac{r^2N_{\nu}^{\prime \prime }(r)+\left(2-2\nu\right)rN_{\nu}^{\prime
				}(r)+\left(\nu^2-\nu\right)N_{\nu}(r)}{rN_{\nu}^{\prime
				}(r)+\left(1-\nu\right)N_{\nu}(r)}=\left(\alpha-1\right)\cos\gamma.
			\end{equation*}
			
			\item[c)] If $\nu \geq \max\{0,\nu_{0} \},$ then the radius of the convex $\gamma-$Spirallike of order $\alpha$ of the function $%
			h_\nu$ is the smallest positive root of the equation 
			\begin{equation*}
			\frac{rN_{\nu}^{\prime \prime }(\sqrt{r})+\left(3-2\nu\right)\sqrt{r}%
				N_{\nu}^{\prime }(\sqrt{r})+\left(\nu^2-2\nu\right)N_{\nu}(\sqrt{r})}{2\sqrt{r%
				}N_{\nu}^{\prime }(\sqrt{r})+2\left(2-\nu\right)N_{\nu}(\sqrt{r})}=\left(\alpha-1\right)\cos\gamma.
			\end{equation*}
		\end{enumerate}
	\end{theorem}	
	\begin{proof}
		We first prove the part $(a)$. From (\ref{eq13}), we have 
		\begin{equation*}
		1+\frac{zf_{\nu }^{\prime \prime }(z)}{f_{\nu }^{\prime }(z)}=1+\frac{
			zN_{\nu }^{\prime\prime}(z)}{N_{\nu }^{\prime}(z)}+\left( \frac{1}{\nu }%
		-1\right) \frac{zN_{\nu }^{\prime}(z)}{N_{\nu }(z)}
		\end{equation*}
		and by means of (\ref{NVZCrp}) and (\ref{NVZTCrp}), we obtain
		\begin{equation*}
		1+\frac{zf_{\nu }^{\prime \prime }(z)}{f_{\nu }^{\prime }(z)}=1-\left( \frac{%
			1}{\nu }-1\right) \sum\limits_{n\geq 1}\frac{2z^{2}}{ \lambda_{\nu ,n}
			^{2}-z^{2}}-\sum\limits_{n\geq 1}\frac{2z^{2}}{\lambda_{\nu ,n}
			^{\prime2}-z^{2}}.
		\end{equation*}
		For $1\geq \nu>\max\left\{0,\nu_0 \right\}$ by using (\ref{eq23}), we get 
		\begin{align}
		& \operatorname{Re}\left(  e^{-i\gamma}\left(1+\frac{zf''_{\nu}(z)}{f'_{\nu}(z)}\right)\right) \nonumber\\
		&= \operatorname{Re}(e^{-i\gamma}) -\operatorname{Re}\left(e^{-i\gamma}\left(\sum_{n\geq1}\frac{2z^2}{{\lambda}^{\prime2}_{\nu, n} -z^2} +\left(\frac{1}{\nu}-1\right)\sum_{n\geq1}\frac{2z^2}{	{\lambda}^2_{\nu, n}-z^2} \right) \right)  \nonumber\\
		&\geq \cos\gamma-\left( \frac{1}{\nu }-1\right) \sum\limits_{n\geq 1}\frac{2r^{2} }{
			\lambda_{\nu,n } ^{2}-r^{2}}-\sum\limits_{n\geq 1}\frac{2r^{2}}{\lambda_{\nu,n} ^{\prime2}-r^{2}} \nonumber\\
		&\geq  \cos\gamma+\frac{r f''_{\nu}(r)}{f'_{\nu}(r)} \label{eq26}
		\end{align}
		where $\left|z\right|=r$.  Moreover, observe that if we use the
		inequality \cite[Lemma 2.1] {Ba1}
		\begin{equation*}
		\mu \operatorname{Re}\left( \frac{z}{a-z}\right) -\operatorname{Re}\left( \frac{z}{b-z}%
		\right) \geq \mu \frac{\left\vert z\right\vert }{a-\left\vert z\right\vert }-%
		\frac{\left\vert z\right\vert }{b-\left\vert z\right\vert }
		\end{equation*}%
		where $a>b>0,$ $\mu \in \lbrack 0,1]$ and $z\in \mathbb{C}
		$ such that $|z|<b$, then we get that the inequality (\ref{eq26}) is also valid
		when $\nu\geq 1$. Here we used that the zeros of $N_{\nu}$ and  $%
		N_{\nu}^\prime$ are interlacing according to Lemma \ref{Le4}. The above
		inequality implies for $r\in (0,\lambda_{\nu ,1}^{\prime})$%
		\begin{equation*}
		\underset{z\in \mathbb{D}_{r}}{\inf }\left\{ \operatorname{Re}\left( e^{-i\gamma}\left(1+\frac{zf''_{\nu}(z)}{f'_{\nu}(z)}\right)-\alpha\cos\gamma\right) \right\} =\left(1-\alpha\right)\cos\gamma+\frac{rf_{\nu }^{\prime \prime }(r)}{f_{\nu }^{\prime }(r)}.
		\end{equation*}%
		Now, the proof of part $(a)$ follows on similar lines as of Theorem \ref{T1}.
		
		For the other parts, note that the functions $g_{\nu}$ and $h_{\nu}$ belong to the Laguerre-P\'{o}lya class $\mathcal{LP}$, which is closed under differentiation, their derivatives $g'_{\nu}$ and $h'_{\nu}$ also belong to $\mathcal{LP}$ and the zeros are real. Thus assuming $\delta_{\nu,n}$  and $\gamma_{\nu,n}$ are the positive zeros of $g'_{\nu}$ and $h'_{\nu}$, respectively, we have the following representations:
		\begin{align*}
		g'_{\nu}(z)=\prod_{n\geq1}\left(1-\frac{z^2}{{\delta}^2_{\nu,n}}\right) \quad
		\text{and} \quad
		h'_{\nu}(z)=\prod_{n\geq1}\left(1-\frac{z}{\gamma_{\nu,n}}\right),
		\end{align*}
		which yield
		\begin{align*}
		1+\frac{zg''_{\nu}(z)}{g'_{\nu}(z)}=1-\sum_{n\geq1}\frac{2z^2}{{\delta}^2_{\nu,n}-z^2} \quad
		\text{and}\quad
		1+\frac{zh''_{\nu}(z)}{h'_{\nu}(z)}=1-\sum_{n\geq1}\frac{z}{{\gamma}_{\nu,n}-z}.
		\end{align*}
		Further, reasoning along the same lines as in Theorem \ref{T1}, the result follows at once. \qed
	\end{proof}
\begin{remark}	
	In Theorem~\ref{T2}, the choice of $\gamma=0$ yields the result for the class of convex functions of order $\alpha$.
\end{remark}

	\begin{table} 
		\caption{Radii of the convex $\frac{\pi}{3}-$Spirallike of order $1/2$ for $f_{1/2 },~g_{1/2}$ and $h_{1/2 }$}
		\begin{tabular}{|c|c|c|c|c|c|c|c|c|c|}
			\hline
			& \multicolumn{3}{c|}{$b=3$ and $c=0$} & \multicolumn{3}{c|}{$a=1$ and $c=0$}
			& \multicolumn{3}{c|}{$a=1$ and $b=2$} \\ \cline{2-10}
			& $a=2$ & $a=3$ & $a=4$ & $b=2$ & $b=3$ & $b=4$ & $c=2$ & $c=3$ & $c=4$ \\ 
			\hline
			\multicolumn{1}{|c|}{  $R^c\left(f_{\frac{1}{2}}\right) $} & $0.0875$ & $0.0677$ & $0.0504$ & 
			$0.0993$ & $0.1131$ & $0.1210$ & 
			$0.1722$ & $0.1925$ & $0.2082$
			\\ \hline
			\multicolumn{1}{|c|}{$R^c\left(g_{\frac{1}{2}}\right)  $} & $0.1221$ & $0.0944$& $0.0703$ 
			& $0.1386$ & $0.1579$ & $0.1691$ & 
			$0.2411$ & $0.2699$ & $0.2921$
			\\ \hline
			\multicolumn{1}{|c|}{$R^c\left(h_{\frac{1}{2}}\right) $} & $0.0406$ & $0.0242$ & $0.0134$ 
			& $0.0524$ & $0.0682$ & $0.0783$ & 
			$0.1622$ & $0.2046$ & $0.2408$
			\\ \hline
		\end{tabular}
	\label{Table2.} 
	\end{table}

	\subsection{Radii of Ma-Minda Starlikeness and Convexity of the Functions $f_{ \nu },$ $g_{\nu }$ and $h_{\nu }$} 	
	In the following, now we establish the result for the complete family of Ma-Minda class of starlike functions. 
	\begin{theorem}\label{Starlike-result}
		Let $\nu \geq \max\{0,\nu_0\}.$ The $\mathcal{S}^*(\varphi)$-radii  of the functions $f_{\nu}$, $g_{\nu}$ and $h_{\nu}$ are respectively, given by the smallest positive root of the equations
		\begin{enumerate}
			\item [$a)$]   $rN'_{\nu}(r)-(1-\beta)\nu N_{\nu}(r)=0$ with $\nu\neq0;$
			\item [$b)$] $ rN'_{\nu}(r)-(\nu-\beta)N_{\nu}(r)=0 ;$
			\item [$c)$] $ \sqrt{r}N'_{\nu}(\sqrt{r})-(\nu-2\beta) N_{\nu}(\sqrt{r})=0$
		\end{enumerate}
		in $|z|<(0, \lambda_{\nu,1})$, $(0,\lambda_{\nu,1})$ and $(0,\lambda_{\nu,1}^{2}),$ where $\varphi(-1)=1-\beta$ and $\beta$ is the radius of the largest disk $\{w : |w-1|<\beta \} \subseteq \varphi(\mathbb{D})$.
	\end{theorem}
	\begin{proof}
		We first prove part $a):$ from equation \eqref{eq22}, we have
		\begin{equation}\label{starlike-expression-f}
		\frac{zf_{\nu }^{\prime }(z)}{f_{\nu }(z)} =\frac{1}{\nu }\frac{ zN_{\nu
			}^\prime(z)}{N_{\nu }(z)}=1-\frac{1}{\nu}\sum\limits_{n\geq 1}\frac{2z^{2}}{
			\lambda_{\nu ,n} ^{2}-z^{2} },\text{ \ \ }
		\left(\nu \geq \max\{0,\nu_0\},~\nu\neq0\right).
		\end{equation}
		Let us consider the continuous function $T_{f_{\nu}}: (0,\lambda_{\nu,1})\rightarrow \mathbb{R}$
		\begin{equation}\label{T-f}
		T_{f_{\nu}}(r)= \frac{1}{\nu}\sum\limits_{n\geq 1}\frac{2r^{2}}{\lambda_{\nu ,n} ^{2}-r^{2} }-\beta.
		\end{equation}
		Then
		\begin{equation*}
		T'_{f_{\nu}}(r)= \frac{1}{\nu}\sum\limits_{n\geq 1} \left(\frac{4r \lambda_{\nu,n}^{2}}{(\lambda_{\nu,n}^{2} -r^2)^2} \right) > 0
		\end{equation*}
		for all $\nu \geq \max\{0,\nu_{0}\}\neq0$ and for $r < \lambda_{\nu,1}$. Also, $T_{f_{\nu}}(0)=-\beta<0$ and 
		$\lim_{r\nearrow \lambda_{\nu ,1}}T_{f_{\nu}}(r)=\infty$. This implies that there exist a unique positive root, say $R_{\varphi}(f_{\nu})$ of the equation $T_{f_{\nu}}(r)=0$ in $(0,\lambda_{\nu,1})$. 	
		Now let $\{w : |w-1|<\beta \} \subseteq \varphi(\mathbb{D})$ such that $\varphi(-1)=1-\beta$. Therefore, using the Lemma~\cite{De} that if $z\in \mathbb{C}$ and $\lambda \in \mathbb{R}$ are such that $ \left\vert z\right\vert\leq r<\lambda,$ then
		\begin{equation*}
		\operatorname{Re}\left( \frac{z}{\lambda -z}\right) \leq \left\vert \frac{z}{\lambda -z}\right\vert  \leq \frac{|z|}{\lambda -\left\vert z\right\vert },
		\end{equation*}
		in \eqref{starlike-expression-f}, we get that
		\begin{align}\label{diskequation-f}
		\left| \frac{zf_{\nu }^{\prime }(z)}{f_{\nu }(z)} -1 \right| =\left| \frac{1}{\nu}\sum\limits_{n\geq 1}\frac{2z^{2}}{\lambda_{\nu ,n} ^{2}-z^{2} }  \right|
		\leq \frac{1}{\nu}\sum\limits_{n\geq 1}\frac{2r^{2}}{\lambda_{\nu ,n} ^{2}-r^{2} }
		\leq \beta
		\end{align}
		implies that $f_{\nu} \in \mathcal{S}^*(\varphi)$ in $|z|< R_{\varphi}(f_{\nu})$. Further, taking $z\in \mathbb{D}$ such that $z=r=-R_{\varphi}(f_{\nu})$, from \eqref{T-f} and \eqref{diskequation-f}, it follows that
		$$	\left| \frac{zf_{\nu }^{\prime }(z)}{f_{\nu }(z)} -1 \right|= \beta,$$
		which implies that for $zf'_{\nu}(z)/f_{\nu}(z) \not\in \varphi(\mathbb{D})$ for all $|z|=r\geq R_{\varphi}(f_{\nu})$, which proves the sharpness part. Reasoning along the same lines, proofs of the parts $b)$ and $c)$ follow. \qed
	\end{proof}	
	
\begin{figure}[h]
		\begin{tabular}{c}
			\includegraphics[scale=0.8]{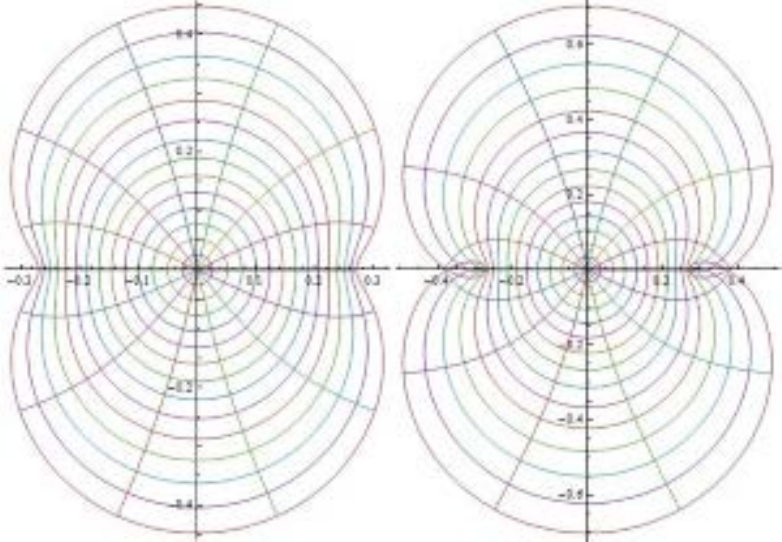}
		\end{tabular}
		\caption{Images of function $g_{1/2}(z)$ for $r=0.3571$ and $r=0.48$, respectively as per Theorem~\ref{Starlike-result}.}\label{Figure2.}
	\end{figure}

		The figure~\ref{Figure2.} was created by taking $a=1,~b=2,~c=0$ and $\varphi(z)=e^z$ in the Theorem~\ref{Starlike-result}.

	\begin{remark}
		Theorem~\ref{Starlike-result} with $\varphi(z)=(1+(1-2\beta)z)/(1-z)$ coincide with Theorem~\ref{T1} with $\alpha=\beta$ and $\gamma=0$.	In particular, choosing $\varphi(z)=(1+z)/(1-z)$ and $\beta=1$ in Theorem~\ref{Starlike-result} reduces to Theorem~\ref{T1} with $\alpha=0$ and $\gamma=0$.
	\end{remark}
	
	As a consequence of Theorem~\ref{Starlike-result}, we get the sufficient conditions for the functions to belong in the class $\mathcal{S}^*(\varphi)$.
	\begin{corollary}
		Let $\nu \geq \max\{0,\nu_0\}.$ Further, $\varphi(-1)=1-\beta$ and $\beta$ is the radius of the largest disk $\{w : |w-1|<\beta \} \subseteq \varphi(\mathbb{D})$. Then
		\begin{enumerate}
			\item [$(i)$] if $1-\frac{1}{\nu}\frac{N'_{\nu}(1)}{N_{\nu}(1)}< \beta$ and $\nu\neq0$, then the function $f_{\nu}$ belongs to $\mathcal{S}^*(\varphi)$.
			\item [$(ii)$] if $\nu-\frac{N'_{\nu}(1)}{N_{\nu}(1)}< \beta$, then the function $g_{\nu}$ belongs to $\mathcal{S}^*(\varphi)$.
			\item [$(iii)$] if $\frac{\nu}{2}-\frac{N'_{\nu}(1)}{2N_{\nu}(1)}< \beta$, then the function $h_{\nu}$ belongs to $\mathcal{S}^*(\varphi)$.
		\end{enumerate} 
	\end{corollary}
	
	\begin{table} 
		\caption{Radii of 
			$\mathcal{S}^*(e^z)$ for the functions  $f_{1/2 },~g_{1/2}$ and $h_{1/2 }$}
		\begin{tabular}{|c|c|c|c|c|c|c|c|c|c|}
			\hline
			& \multicolumn{3}{c|}{$b=3$ and $c=0$} & \multicolumn{3}{c|}{$a=1$ and $c=0$}
			& \multicolumn{3}{c|}{$a=1$ and $b=2$} \\ \cline{2-10}
			& $a=2$ & $a=3$ & $a=4$ & $b=2$ & $b=3$ & $b=4$ & $c=2$ & $c=3$ & $c=4$ \\ 
			\hline
			\multicolumn{1}{|c|}{  $r^*_{e^z}\left(f_{\frac{1}{2}}\right)  $} & $0.2354$ & $0.1818$
			& $0.1352$ & $0.2674$ & $0.3050$ & $0.3268$ & $0.4696$ & $0.5272$ & $0.5718$
			\\ \hline
			\multicolumn{1}{|c|}{$r^*_{e^z}\left(g_{\frac{1}{2}}\right)  $} & $0.3140$ & $0.2421$& $0.1798$ 
			& $0.3571$ & $0.4078$ & $0.4372$ & 
			$0.6350$ & $0.7160$ & $0.7792$
			\\ \hline
			\multicolumn{1}{|c|}{$r^*_{e^z}\left(h_{\frac{1}{2}}\right)  $} & $0.1611$ & $0.0953$& $0.0523$ 
			& $0.2090$ & $0.2736$ & $0.3150$ 
			& $0.6845$ & $0.8819$ & $1.0559$
			\\ \hline
		\end{tabular}
		\label{Table3.} 
	\end{table}		
	
	In the following, we prove the result for the Ma-Minda class of convex functions.
	\begin{theorem}\label{convex-result}
		The $\mathcal{C}(\varphi)$-radii  of the functions $f_{\nu}$, $g_{\nu}$ and $h_{\nu}$ are respectively, given by the smallest positive root of the equations
		\begin{enumerate}
			\item [$a)$] $rf''_{\nu}(r)+\beta\nu f'_{\nu}(r)=0,$ where $\nu \geq \max\{0,\nu_0\},~\nu\neq0, $
			\item [$b)$] $r g''_{\nu}(r)+\beta g'_{\nu}(r)=0,$ where $\nu \geq \max\{0,\nu_0\},$
			\item [$c)$] $\sqrt{r} h''_{\nu}(\sqrt{r})+\beta h'_{\nu}(\sqrt{r})=0,$ where $\nu \geq \max\{0,\nu_0\}$
		\end{enumerate}
		in $|z|<(0, \lambda_{\nu,1})$, $(0,\lambda_{\nu,1})$ and $(0,\lambda_{\nu,1}^{2}),$ where $\varphi(-1)=1-\beta$ and $\beta$ is the radius of the largest disk $\{w : |w-1|<\beta \} \subseteq \varphi(\mathbb{D})$.
	\end{theorem}
	\begin{proof}
		From \eqref{starlike-expression-f}, it follows that
		\begin{equation}\label{convex-expression-f}
		\frac{zf''_{\nu}(z)}{f'_{\nu}(z)} = \frac{zN''_{\nu}(z)}{N'_{\nu}(z)}+ \left(\frac{1}{\nu}-1 \right) \frac{zN'_{\nu}(z)}{N_{\nu}(z)}.
		\end{equation}
		Now using the Weierstrass decompositions of the functions $N_{\nu}$ and $N'_{\nu}$, which are given below
		\begin{equation*}
		N_{\nu}(z)= \frac{Q(\nu) z^{\nu}}{2^{\nu}\Gamma(\nu+1)} \prod\limits_{n\geq1}\left( 1- \frac{z^2}{\lambda_{\nu,n}^{2}} \right)
		\end{equation*}
		and 
		\begin{equation*}
		N'_{\nu}(z)= \frac{Q(\nu) \nu z^{\nu-1}}{2^{\nu}\Gamma(\nu+1)} \prod\limits_{n\geq1}\left( 1- \frac{z^2}{\tilde{\lambda}_{\nu,n}^{2}} \right),
		\end{equation*}
		where $\lambda_{\nu,n}$ and $\tilde{\lambda}_{\nu,n}$ are the $n^{th}$ positive roots of $N_{\nu}$ and $N'_{\nu}$, we have from \eqref{convex-expression-f} that
		\begin{equation}\label{f-series-expression}
		1+ \frac{zf''_{\nu}(z)}{f'_{\nu}(z)} =1- \left(\frac{1}{\nu}-1 \right) \sum\limits_{n\geq1} \frac{2z^2}{\lambda_{\nu,n}^{2}-z^2} -\sum\limits_{n\geq1}\frac{2z^2}{\tilde{\lambda}_{\nu,n}^{2}-z^2}.
		\end{equation}
		Using the Lemma~\cite{De} in \eqref{f-series-expression} that if $z\in \mathbb{C}$ and $\lambda \in \mathbb{R}$ are such that $ \left\vert z\right\vert\leq r<\lambda,$ then
		\begin{equation*}
		\operatorname{Re}\left( \frac{z}{\lambda -z}\right) \leq \left\vert \frac{z}{\lambda -z}\right\vert  \leq \frac{|z|}{\lambda -\left\vert z\right\vert },
		\end{equation*} 
		it follows that the inequalities
		\begin{equation*}
		\left| \frac{zf''_{\nu}(z)}{f'_{\nu}(z)} \right| \leq \left(\frac{1}{\nu}-1 \right) \sum\limits \frac{2r^2}{\lambda_{\nu,n}^{2}-r^2} + \sum\limits_{n\geq1} \frac{2r^2}{\tilde{\lambda}_{\nu,n}^{2}-r^2} 
		= -\frac{rf''_{\nu}(r)}{f'_{\nu}(r)}
		\end{equation*}
		hold in $|z|=r<\tilde{\lambda}_{\nu,1}$ for $\nu\leq1$. Observe that the inequality $\left|{zf''_{\nu}(z)}/{f'_{\nu}(z)} \right| \leq -{rf''_{\nu}(r)}/{f'_{\nu}(r)}$ also holds for $\nu >1$ in $|z|=r<\tilde{\lambda}_{\nu,1}$. Let $R_{\varphi}^{c}(f_{\nu})$ be the smallest positive root of the equation
		$$\frac{rf''_{\nu}(r)}{f'_{\nu}(r)}+\beta=0.$$
		Now let $\{w : |w-1|<\beta \} \subseteq \varphi(\mathbb{D})$ such that $\varphi(-1)=1-\beta$. Then it follows that 
		\begin{equation*}
		\left| \frac{zf''_{\nu}(z)}{f'_{\nu}(z)} \right| 
		\leq  -\frac{rf''_{\nu}(r)}{f'_{\nu}(r)}
		\leq \beta
		\end{equation*}
		hold in $|z|=r<R_{\varphi}^{c}(f_{f_{\nu}})$ for $\nu \geq \max\{0,\nu_0\},~\nu\neq0$, which implies that the function $f_{\nu} \in \mathcal{C}(\varphi)$ in $|z|=r<R_{\varphi}^{c}(f_{f_{\nu}})$ for $\nu \geq \max\{0,\nu_0\},~\nu\neq0.$ For the sharpness, let us consider $z=-R_{\varphi}^{c}(f_{\nu})$. Then 
		$$\left| \frac{zf''_{\nu}(z)}{f'_{\nu}(z)} \right| 
		=  -\frac{rf''_{\nu}(r)}{f'_{\nu}(r)}
		= \beta$$
		such that  $1+zf''_{\nu}(z)/f'_{\nu}(z) \not\in \varphi(\mathbb{D})$ for all $|z|=r\geq R_{\varphi}^{c}(f_{\nu})$, which proves the sharpness part. Reasoning along the same lines, proofs of the parts $b)$ and $c)$ follow.  	\qed
	\end{proof}	

\begin{table}
	\caption{Radii of 
		$\mathcal{C}(\varphi)$ for the functions $f_{1/2 },~g_{1/2}$ and $h_{1/2 }$ with $\varphi(z)=z+\sqrt{1+z^2}$}
	\begin{tabular}{|c|c|c|c|c|c|c|c|c|c|}
		\hline
		& \multicolumn{3}{c|}{$b=3$ and $c=0$} & \multicolumn{3}{c|}{$a=1$ and $c=0$}
		& \multicolumn{3}{c|}{$a=1$ and $b=2$} \\ \cline{2-10}
		& $a=2$ & $a=3$ & $a=4$ & $b=2$ & $b=3$ & $b=4$ & $c=2$ & $c=3$ & $c=4$ \\ 
		\hline
		\multicolumn{1}{|c|}{${r}^c_{\varphi}(f_{\frac{1}{2} })$} & $0.1271$ & $0.0983$& $0.0732$ 
		& $0.1443$ & $0.1644$ & $0.1761$ 
		& $0.2512$ & $0.2812$ & $0.3044$
		\\ \hline
		\multicolumn{1}{|c|}{${r}^c_{\varphi}(g_{\frac{1}{2}  })$} & $0.1749$ & $0.1352$& $0.1005$ &
		$0.1987$ & $0.2266$ & $0.2427$ &
		$0.3482$ & $0.3906$ & $0.4234$
		\\ \hline
		\multicolumn{1}{|c|}{${r}^c_{\varphi}(h_{\frac{1}{2}  })$} & $0.0759$ & $0.0451$ & $0.0248$ 
		& $0.0983$ & $0.1283$ & $0.1474$ 
		& $0.3124$ & $0.3979$ & $0.4719$
		\\ \hline
	\end{tabular}
	\label{Table4.} 
\end{table}

The following corollary yields sufficient conditions on parameters for functions to be in $\mathcal{C}(\varphi)$.
	\begin{corollary}
		Let $\nu \geq \max\{0,\nu_0\},~\nu\neq0 $. Let $\varphi(-1)=1-\beta$ and $\beta$ is the radius of the largest disk $\{w : |w-1|<\beta \} \subseteq \varphi(\mathbb{D})$. Then
		\begin{enumerate}
			\item [$(i)$] The function $f_{\nu}$ belongs to $\mathcal{C}(\varphi)$, if 
			$$\left(1-\frac{1}{\nu}\right)\frac{N'_{\nu}(1)}{N_{\nu}(1)} -\frac{N''_{\nu}(1)}{N'_{\nu}(1)} < \beta \quad \text{and} \quad 1\geq\nu \geq \max\{0,\nu_0\},~\nu\neq0 \quad \text{and}\quad \nu>1.$$ 
			\item [$(ii)$] The function $g_{\nu}$ belongs to $\mathcal{C}(\varphi),$ if 
			$$\frac{N''_{\nu}(1)+2(1-\nu)N'_{\nu}(1)+({\nu}^2+\nu)N_{\nu}(1)}{(1-\nu)N_{\nu}(1)-N'_{\nu}(1)} < \beta$$
			\item [$(iii)$] The function $h_{\nu}$ belongs to $\mathcal{C}(\varphi)$, if 
			$$ \frac{1}{2}\frac{N''_{\nu}(1)+ (3-2\nu)N'_{\nu}(1)+ ({\nu}^2-2\nu)N_{\nu}(1)}{(\nu-2)N_{\nu}(1)-N'_{\nu}(1)} < \beta. $$
		\end{enumerate} 
	\end{corollary}

\begin{figure}[h]
	\begin{tabular}{c}
		\includegraphics[scale=0.9]{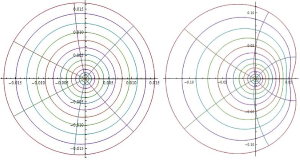}
	\end{tabular}
	\caption{Images of function $h_{1/2}(z)$ for $r=0.1271$ and $r=0.3$, respectively as per Theorem~\ref{convex-result} and Remark~\ref{manyclass}.}\label{Figure3.}
\end{figure}

The figure~\ref{Figure3.} was created by taking $\nu=1/2,~a=2,~b=3,~c=0$ and $\varphi(z)=z+\sqrt{1+z^2}$ in Theorem~\ref{convex-result}.

	\begin{remark}\label{manyclass}
		Let $\beta$ be the radius of the lagest disk $\{w: |w-1|< \beta\}$ inside $\varphi(\mathbb{D})$. Then Theorem~\ref{Starlike-result} and Theorem~\ref{convex-result} hold for each of the following cases and the radii are sharp, where
		\begin{enumerate}
			\item [$(i)$]
			$\beta=\min\left\{\left|1-\frac{1+A}{1+B}\right|, \left|1-\frac{1-A}{1-B}\right|\right\}=\frac{A-B}{1+|B|}$ when $\varphi(z)= \frac{1+Az}{1+Bz}$, where $-1\leq B<A\leq1$;
			
			\item  [$(ii)$]
			$\beta=\sqrt{2-2\sqrt{2}+\sqrt{-2+2\sqrt{2}}}$ when  $\varphi(z)=\sqrt{2}-(\sqrt{2}-1)\sqrt{\frac{1-z}{1+2(\sqrt{2}-1)z}}$;
			
			\item  [$(iii)$]
			$\beta=\sqrt{2}-1$ when $\varphi(z)=\sqrt{1+z}$;
			
			\item [$(iv)$]
			$\beta=1-{1}/{e}$ when $\varphi(z)=e^z$;
			
			\item [$(v)$]
			$\beta=2-\sqrt{2}$ when $\varphi(z)=z+\sqrt{1+z^2}$;
			
			\item [$(vi)$]
			$\beta=\frac{e-1}{e+1}$ when $\varphi(z)=\frac{2}{1+e^{-z}}$;
			
			\item  [$(vii)$]
			$\beta=\sin{1}$ when $\varphi(z)=1+\sin{z}$;
			
			\item [$(viii)$]
			$\beta=1-e^{e^{-1}-1}$ when $\varphi(z)=e^{e^z -1}$;
			
			\item [$(ix)$]
			for the domains bounded by the conic sections
			$\Omega_\kappa:=\{w=u+iv: u^2>{\kappa}^2(u-1)^2+{\kappa}^2v^2; \kappa\in[0,\infty)  \},$ we have $$\beta=\frac{1}{\kappa+1},$$ 
			where the boundary curve of $\Omega_\kappa$ for fixed $\kappa$ is represented by the imaginary axis $(\kappa=0)$, the right branch of a hyperbola $(0<\kappa<1)$, a parabola $(\kappa=1)$ and an ellipse $(\kappa>1)$. The univalent Carath\'{e}odory functions mapping $\mathbb{D}$ onto $\Omega_\kappa$ is given by
			\begin{equation*}
			\varphi(z):=\varphi_{\kappa}(z)= 
			\left\{
			\begin{array}{lll}
			\frac{1+z}{1-z} & $for$ & \kappa=0;\\
			1+\frac{2}{1-\kappa^2}\sinh^2(A(\kappa) arctanh\sqrt{z}) & $for$ & \kappa\in(0,1);\\
			1+\frac{2}{\pi^2}\log^2{\frac{1+\sqrt{z}}{1-\sqrt{z}}} & $for$ & \kappa=1;\\
			1+\frac{2}{\kappa^2-1}\sin^2\left(\frac{\pi}{2K(t)}F\left( \frac{\sqrt{z}}{\sqrt{t}}, t \right) \right) & $for$ & \kappa>1,		
			\end{array}	
			\right.
			\end{equation*}
			where $A(\kappa)=(2/\pi)\arccos(\kappa)$, $F(w,t)=\int_{0}^{w}\frac{dx}{\sqrt{(1-x^2)(1-t^2x^2)}}$ is the Legender elliptic integral of the first kind, $K(t)=F(1,t)$ and $t\in(0,1)$ is choosen such that $\kappa=\cosh(\pi K'(t)/2K(t))$.
		\end{enumerate}
	\end{remark}

\section*{Statements and Declarations}
\begin{itemize}
	\item {\bf Funding}: The work of Kamajeet Gangania is supported by University Grant Commission, New-Delhi, India  under UGC-Ref. No.:1051/(CSIR-UGC NET JUNE 2017).
	\item {\bf Conflict of interest}: The authors declare that they have no conflict of interest
	\item {\bf Availability of data and materials }: NA
	\item {\bf Code availability}: NA 
	\item {\bf Authors' contributions }: All authors contributed Equally.
\end{itemize}

\end{document}